\documentclass[]{amsart}

\usepackage{cite}
\usepackage{a4wide}
\usepackage{amssymb}

\newtheorem{theorem}{Theorem}[section]
\newtheorem{lemma}[theorem]{Lemma}
\newtheorem{proposition}[theorem]{Proposition}
\newtheorem{corollary}[theorem]{Corollary}

\theoremstyle{definition}
\newtheorem{definition}[theorem]{Definition}

\theoremstyle{remark}
\newtheorem{remark}[theorem]{Remark}

\newcommand\grad[1]{\nabla #1}
\newcommand\norma[1]{\left\|#1\right\|}
\newcommand\nq[1]{\left\|#1\right\|_{\Q}}
\newcommand\nqmu[1]{\left\|#1\right\|_{\Q_\mu}}
\newcommand\nll[1]{\left\|#1\right\|_{\ell^2}}

\newcommand\nlmu[1]{\left\|#1\right\|_{\ell^2_\mu}}

\def\Z{\mathbb{Z}}
\def\R{\mathbb{R}}

\def\Id{\mathbb I}

\def\C{{\mathcal C}}

\def\E{{\mathcal E}}
\def\U{{\mathcal U}}
\def\L{{\mathcal L}}
\def\I{{\mathcal I}}
\def\Q{{\mathcal Q}}

\def\tond#1{\left(#1\right)}
\def\quadr#1{\left[#1\right]}
\def\inter#1{\langle#1\rangle}

\begin{document}

\bibliographystyle{siam}

\title[Small Amplitude breathers in KG lattices]
{Existence and continuous approximation\\of small amplitude
  breathers\\in 1D and 2D Klein--Gordon lattices}

\author[D.~Bambusi]{Dario Bambusi}
\author[S.~Paleari]{Simone Paleari}
\author[T.~Penati]{Tiziano Penati}
  
\email{dario.bambusi@unimi.it}
\email{simone.paleari@unimi.it}
\email{tiziano.penati@unimi.it}

\address{Dipartimento di Matematica ``F. Enriques'', Via Saldini 50,
  20133 Milano, Italy.}

\subjclass[2000]{34C25, 37K50, 37K60, 65N30.}

\keywords{Breathers, small amplitude solutions, 1D and 2D,
  Klein--Gordon lattices, approximation and interpolation.}

\begin{abstract}
We construct small amplitude breathers in 1D and 2D Klein--Gordon
infinite lattices. We also show that the breathers are well
approximated by the ground state of the nonlinear Schr\"odinger
equation. The result is obtained by exploiting the relation between
the Klein Gordon lattice and the discrete Non Linear Schr\"odinger
lattice. The proof is based on a Lyapunov-Schmidt decomposition and
continuum approximation techniques introduced in \cite{BamP09},
actually using its main result as an important lemma.
\end{abstract}

\maketitle


\section{Introduction}
\label{s.i}

In the present paper we prove existence of small amplitude breathers
in some nonlinear Klein Gordon Lattices in dimension one and
two. Moreover we prove that such solutions are well approximated by
the ground state of a suitable nonlinear Schr\"odinger
equation.

The present paper is a direct continuation of \cite{BamP09} where the
same result was obtained for the discrete nonlinear Schr\"odinger
equation. More precisely, in \cite{BamP09} two of us proved the
existence of spatially localized, time periodic solutions in one and
two dimensional discrete Nonlinear Schr\"odinger equation (dNLS). In
particular the so called Sievers-Takeno (ST) and Page (P) modes in 1D,
and also the hybrid (H) modes in 2D were constructed. The breathers of
\cite{BamP09} have been obtained as critical points of the Energy
functional constrained to the surface of constant $\ell^2$ norm. In
turn they were constructed by continuation from the ground state of
the Nonlinear Schr\"odinger equations (NLS).  Thus such breathers
turned out to be well approximated by the corresponding solutions of
the continuous model. A key role in the proof was played by the
interpolation of sequences (configurations of the discrete system) by
the use of the so called Finite Elements, usually used in numerical
analysis.

In order to apply those ideas to the Klein Gordon lattice in which the
trivial variational characterization of the breathers is no more true
(see however \cite{MW89} for general results on periodic orbits via
variational methods, and \cite{Bak04,BakP04} for recent results on KG
chains), we have to establish a new connection between the KG lattice
and the dNLS lattice. Actually the dNLS is heuristically known to be a
resonant normal form (or modulation equation) of the KG lattice, an
idea which has been exploited in order to describe the finite time
dynamics of KG lattices. However, we need here a connection suitable
for the description of infinite time dynamics. Such a connection is
obtained by using the method of Lyapunov-Schmidt decomposition, in
which the so called Kernel equation turns out to be a perturbation of
the dNLS.  We recall that a connection between normal form theory and
the method of Lyapunov-Schmidt decomposition was first recognized in
\cite{Dui84} and exploited in \cite{BP,BP2}.

From a technical point of view the procedure is quite delicate, and in
order to obtain a meaningful result we have to exploit the
techniques introduced in \cite{BamP09}. We recall that in such a
paper the idea was to use the Finite Elements in order to interpolate
a sequence (configuration of the lattice) with a function of class
$H^1$. This allows to consider functionals on the discrete
configuration space as restrictions to suitable subspaces of
functionals on the continuous phase space.

This paper is part of a research line consisting in studying the
dynamics of lattices using the continuous approximation. Along this
line many results on the finite time dynamics have been proved
\cite{Kal89,KirSM92,SchW00,BamCP02,BPon05,BamCP09}, but little is
known on the approximation for infinite times (see
\cite{FriP99,FriP1-4,HofW08,MizP08} on the dynamics close to solitons
in FPU type models, and the papers
\cite{Ioo00,IooK00,IooJ05,IooP06,Jam03,JamS05,JamS08,Sir05,BamP09} for
what concerns existence of breathers).

We recall that existence of breathers in lattices has been proved in
\cite{MacA94} (see also \cite{Bam96,SepM97}), by looking at the
opposite limit, namely the anti-continuous one which leads to large
amplitude breathers. Concerning small amplitude breathers, their
existence have been established in one dimensional lattices using the
spatial dynamics approach (see, e.g.,
\cite{Ioo00,IooK00,IooJ05,IooP06,Jam03,JamS05,JamS08,Sir05}). In the
case of higher dimensional lattices we recall the result by Aubry,
Kopidakis, Kadelburg \cite{AubKK01} (which in principle should be
applicable also to the present model) and the results by Weinstein
\cite{Wei99} (not directly applicable to the present model). However
all these results are obtained by topological methods and give no
information on the number of existing breathers and on their
shape. The result of the present paper is actually the first one
allowing to explicitly construct the ST, as well as the P and the H
modes.

\section{Main result}
\label{ss.mr}

We consider the equation:
\begin{equation}
\label{e.KG}
\ddot{q_j} = a(\Delta q)_j - V'(q_j) \ ,\quad j\in\mathbb{Z}^n
\end{equation}
with
\[
  V(s)=\frac{s^2}2 - \frac{\beta}{2p+2}|s|^{2p+2} \, ,\quad   (\Delta
  q)_j := \sum_{k:|k-j|=1}(q_k-q_j) \, , 
\]
with $|k|=\sum_{m=1}^n|k_m|$.

Here, and in the following, we will always consider $n\in \{1,2\}$.

\begin{remark}
The value of the constant $\beta$ is not essential since it can be
changed by means of an amplitude rescaling. For our convenience we fix
it as $\beta=c_1^{-1}$, with $c_1$ defined in
\eqref{e.c1}.\end{remark}

\begin{remark}
Due to the choice of the sign in front of the nonlinearity we obtain
that the continuous approximation of the normal form is the {\it
  focusing} NLS equation. This is crucial for our analysis which does
not apply to the defocusing case.
\end{remark}

To state the approximation part of our result, we need to refer to the
ground state of the Nonlinear Schr\"odinger equation. Consider the
Nonlinear Schr\"odinger Equation (NLS) $$\imath
\dot\psi=-\Delta_c\psi-|\psi|^{2p}\psi$$ in $\R^n$,
where $$\Delta_c:=\sum_{j=1}^n \frac{\partial^2}{\partial x_j^2}$$ is
the usual Laplacian operator. It is well known that, if $p<2/n$ there
exists a unique ground state $\psi_c$ of the NLS fulfilling the
additional requirements of being real valued, positive, radially
symmetric and exponentially decaying
\cite{BerL83,ColGM78,GriSS87}. Such a ground state is defined as the
function which realizes the minimum of
$\int_{\R^n}\left[|\nabla\psi|^2-\frac1{p+1}|\psi|^{2p+2}\right]$
restricted to $\int_{\R^n}|\psi|^2=1$. It can often be computed or
described quite explicitly.

For any $\mu>0$ small enough consider the following $2^n$ distinct
sequences $\psi^i\equiv\left\{\psi_j^i(\mu)\right\}_{j\in\Z^n}$ defined
by restricting the NLS ground state $\psi_c$ onto $\Z^n$ as follows

\begin{equation}
\label{e.Psiref}
\begin{cases} \psi^1_{j}:=
\psi_c\tond{\mu j},\\ \psi^2_{j}:=\psi_c\tond{\mu j +
  \frac\mu2},
\end{cases}  n=1;\qquad\qquad
\begin{cases} \psi^1_{j_1,j_2}:=\psi_c\tond{\mu j_1, \mu
      j_2},\\ \psi^2_{j_1,j_2}:=\psi_c\tond{\mu j_1, \mu j_2 +
      \frac\mu2},\\ \psi^3_{j_1,j_2}:=\psi_c\tond{\mu j_1 + \frac\mu2, \mu
      j_2},\\ \psi^4_{j_1,j_2}:=\psi_c\tond{\mu j_1 + \frac\mu2, \mu j_2 +
      \frac\mu2},
\end{cases}  n=2.
\end{equation}
These reference sequences correspond to the ST and P modes, plus the H
modes in two dimensions.

As a last step, for all the previous sequences, we renormalize the
amplitude and add a temporal dependence in the following way:
\begin{equation}
\label{e.refsol}
  \Psi^i(t)=\mu^{\frac1p}\cos(\omega t)\psi^i \, ,
\end{equation}
with $\omega=\omega(\mu):=\sqrt{1-m\mu^2}$, where $m$ is a real
constant (see sect. \ref{ss.LS}).

\begin{remark}
\label{r.profil}
By construction the sequences $\psi^i$ are uniformly bounded in $\mu$,
but with diverging $\ell^2$ norm as $\mu\to 0$;
so we are calling ``breathers'' solutions which are localized on an
increasing interval $[-k,k]^n$ with $k\sim 1/\mu$. The reference
solution \eqref{e.refsol} share the same localization property, but
with bounded $\ell^2$ norm, due to the amplitude rescaling.

\end{remark}

We are now ready to state our result.

\begin{theorem}
\label{t.main}
Assume $n\in\{1,2\}$, $0<a<\frac12$ and $\frac12\leq p<\frac{2}{n}$,
then there exists $\mu^*>0$, such that for any $0<\mu<\mu^*$ there
exist $2^n$ distinct real valued sequences
$q^i\equiv\left\{q^i_j\right\}_{j\in \Z^n}\in H^2([0,T];\ell^2)$,
which are time periodic solutions of \eqref{e.KG} with period
$T=\frac{2\pi}{\omega(\mu)}$. Such solutions fulfill
\begin{eqnarray}
\label{e.mest}
  \norma{q^i - \Psi^i}_{H^2([0,T];\ell^2)} &\leq C_2\mu^r \,,\qquad
  r&:=\frac{1}{p}-\frac{n}2+1 \; ;
\\
\label{e.supest}
\sup_{j}  \left|{q^i_j(t) - \Psi^i_j(t)} \right| &\leq C_2\mu^s \,,\qquad
  s&:=\frac{1}{p}-\frac{n}2+\frac32 \; .
\end{eqnarray}
\end{theorem}

\begin{remark}The periodic orbits we find are actually $\C^{3,1}$ in time,
  thus they are classical solutions. Indeed, once we get they are
  $H^2$, and consequently $\C^{1,1}$ by Sobolev embeddings, since the
  operator $q\to |q|^{2p}q$ maps $\C^{1,1}$ into itself, one has
  $\ddot q\in\C^{1,1}$ from equation \eqref{e.KG}.
\end{remark}

\begin{remark}
 The bound $p<2/n$ comes from an analogous bound for the existence and
 approximability of the ground state of the dNLS (see \cite{BamP09}).
 On the other hand, for the dNLS, Weinstein proved in \cite{Wei99} the
 non existence of small amplitude ground states for $p\geq\frac2n$.
\end{remark}

\begin{remark}The reference solution, when measured in
  ${H^2([0,T];\ell^2)}$ has norm of order
  $\mu^{\frac{1}{p}-\frac{n}{2}}$, so the estimate \eqref{e.mest}
  shows that in such a norm the distance between the actual solution
  and the reference solution is small compared to the size of the
  solution. On the contrary such an estimate gives no information on
  the distance between the single particle in the approximate and the
  true solution (sup norm). A better and relevant control is given by
  the estimate \eqref{e.supest}, which is obtained through the use of
  a discrete analogue of Sobolev embedding theorems (see
  Sect. \ref{s.eop}).
\end{remark}

\begin{remark}We also stress that our approximation estimates, not only
  control the spatial profile, but also the time dependence. In
  particular, one easily gets from \eqref{e.wnv} and \eqref{nv.1} that
  the first harmonic gives the principal contribution, the others
  being small corrections of order $\mu^\sigma$, with
  $\sigma:=\frac1p-\frac{n}2+2$.
\end{remark}

The rest of the paper is devoted to the proof of theorem \ref{t.main},
and is organized as follows. Section~\ref{s.s} contains the setting of
the problem, with the description of the proof of the main
Theorem. The range equation is dealt with in Section~\ref{s.R}, while
the solution of the kernel one is discussed in Section~\ref{s.Q};
final estimates are presented in Section~\ref{s.eop}.  Some technical
details are given in Appendixes: in Appendix~\ref{s.nem} we show some
regularity results for the nonlinearity; in Appendix~\ref{s.A.norms}
we give some improved estimates on the approximation of $\ell^q$
norms; and in Appendix~\ref{app.3} we prove the extension of the
Implicit Function Theorem used to prove Proposition~\ref{t.K}.


\section{Settings and proof of Theorem \ref{t.main}}
\label{s.s}

In all the paper we will deal only with sequences $q_j$ which are
reflection invariant, and thus which fulfill $q_{j}=q_{-j}$. Thus,
{\sl when writing $\ell^2$ we will actually mean the subspace of
  $\ell^2$ composed by symmetric sequences. The same will be true for
  all the other spaces of sequences that we will meet in the paper.}

More precisely, we will denote the space by $\ell^2$ whenever it is
endowed with its standard scalar product and norm, namely
\[
  \inter{q,p}_{\ell^2}:=\sum_j q_j p_j \, , \qquad
  \nll{q}^2:=\inter{q,q}_{\ell^2} \, ;
\]
and it will be denoted by $\Q$ when endowed with the
norm
\begin{equation}
\label{e.star}
  \nq{q}^2 := \nll{q}^2 + \frac1{\mu^2}\inter{q,-\Delta q}_{\ell^2}^2 \, ,
\end{equation}
which will play a fundamental role in Lemma~\ref{nv} and
Section~\ref{s.Q}; here $\mu>0$ is the small parameter which was
introduced in \eqref{e.refsol}.

\begin{remark}
From the technical viewpoint, the interplay between the norms
$\nll{\cdot}$ and $\nq{\cdot}$ is one of the delicate points. Indeed
it turns out that the $\Q$ norm is too strong to ensure enough
regularity for all the continuation procedures; but the $\ell^2$ one
is too week to grasp all the relevant information, in particular the
non-degeneracy of the solutions, contained in the result of
\cite{BamP09}. We will thus play with both the norms depending on the
situations.
\end{remark}

We look for $T$-periodic solutions of \eqref{e.KG} 
of the form
\[
  q_j(t)=u_j(\omega t) \, ,
\]
with\footnote{\ Due to the autonomous and reversible nature of
  \eqref{e.KG}, it is rather natural to look for solutions even in
  time, thus with a Fourier development in cosine only.}
\begin{equation}
\label{e.ujlel}
  u_j(t)=\sum_{l\ge0}u_{j,l}\cos(lt) \,.
\end{equation}

We define $X_k\subset H^k([0,2\pi];\ell^2)$, to be the space of {\it even}
functions of time taking value in $\ell^2$ which are square
integrable together with their weak derivatives up to order $k$.

In
this paper we will only use $k\in\{0,2\}$; for $k=2$ we
will use the norm
\begin{equation}
\label{e.norm}
\norma{u}^2_{X_2} := \int_0^{2\pi}\tond{\nll{u(t)}^2 + \nll{\dot u(t)}^2
+ \nll{\ddot u(t)}^2}dt  = \sum_{j,l}(1+l^2+l^4)u_{j,l}^2 \, .
\end{equation}

Using the time-rescaled variable $u$, we rewrite \eqref{e.KG} in the form
\begin{equation}
\label{e.LN}
  L^{(\omega)} u = N(u) \, ,
\end{equation}
where
\begin{eqnarray}
\label{e.LL}
  L^{(\omega)} &:&= L_\omega - a\Delta \, , \qquad\qquad
  L_\omega := \omega^2\partial_{tt} + \Id \, ;
\\
\label{e.NN}
  \left(N(u)\right)_j &:&=  \beta |u_j|^{2p}u_j  \, , 
\end{eqnarray}
are respectively the linear and nonlinear operators.


\subsection{Lyapunov--Schmidt decomposition}
\label{ss.LS}

We look for small amplitude solutions, bifurcating from the trivial
one.  We perform a Lyapunov--Schmidt decomposition with respect to
$L_1\equiv L_\omega\big\vert_{\omega=1}$. Define
\begin{equation}
\label{e.VR}
  V_2 := Ker (L_1) \, , \qquad\quad
  W_0 := Range (L_1) \, ,\qquad W_2:= W_0\cap X_2
\end{equation}
and $V_0$, the completion of $V_2$ in $X_0$.  Such decomposition is
invariant under the action of the linear operators \eqref{e.LL}.

We correspondingly decompose $u$ in the two components $v$ and $w$ 
\begin{equation}
\label{e.vw}
  u=v+w \,, \qquad\quad u\in X_2 \,, \quad v\in V_2 \,, \quad w\in W_2 \,;
\end{equation}
according to the development \eqref{e.ujlel}, if we denote
$e_l=\cos(lt)$ and $\L:=\{ l \in \Z, l\geq 0, l\neq 1 \}$, we
may write
\begin{equation}
\label{e.uvw}
  u_j(t)=v_j e_1 + \sum_{l\in\mathcal{L}}w_{j,l} e_l \,,
\end{equation}
where $v_j=v_{j,1}$ is the only Fourier component in the
kernel. Denote $\lambda := \mu^2m=1-\omega^2$, $\Pi_V$ the projector
onto $V_0$, and $\Pi_W=\Pi_V-\mathbb{I}$ the projector onto $W_0$.  We
remark that $V_2$ and $V_0$ are isometrically (up to a constant factor)
isomorphic to $\ell^2$, since any element $v\in V_2$ or $z\in V_0$
have only one Fourier component.

\subsection{Proof of Theorem \ref{t.main}}

In this subsection we give the steps of the proof, leaving
the full details to the subsequent sections.

\subsubsection*{Step 0: Decomposition}
We first decompose \eqref{e.LN} in the following two equations, the
first on the range and the second on the kernel
\begin{eqnarray}
\label{e.R}
  L^{(\omega)} w &=& \Pi_W(N(v+w)) \, ,
\\
\label{e.K}
  \lambda v - a(\Delta v) &=& \Pi_V(N(v+w))  \, .
\end{eqnarray}

\subsubsection*{Step 1: Range equation}
As usual in the Lyapunov--Schmidt decomposition, we first fix $v$ and
$\omega$ and solve the range equation \eqref{e.R} via the Implicit
Function Theorem (see Proposition~\ref{t.R}) showing that
\[
  w(v) \simeq \left(L^{(\omega)}\right)^{-1}
  \Pi_W(N(v))+\text{higher\ order\ terms} \, .
\]

We stress that the range equation is solved using always the $\ell^2$
norm. All the corresponding discussions and proofs are presented in
Section \ref{s.R}.

\subsubsection*{Step 2: Kernel equation}

We insert now the solution of the range equation $w(v)$ in the kernel
equation \eqref{e.K}.  Due to the smallness of $w(v)$, the term
$\Pi_VN(v+w(v))$ can be split into a main part $\Pi_VN(v)$ and a
remainder:

\begin{equation}
\label{e.q0q1}
  \Pi_VN(v+w(v)) = \Pi_VN(v) + \Pi_V\quadr{N(v+w(v)) - N(v)} \, ;
\end{equation}
it turns out that $\tond{\Pi_VN(v(t))}_j = \beta c_1|v_j|^{2p} v_j
\cos(t)$, and
\begin{equation}
\label{e.c1}
c_1:=\int_0^{2\pi}|\cos(t)|^{2p}\cos^2(t)dt \; .
\end{equation}
By applying the following scaling
\begin{equation}
\label{e.scalK}
\lambda = m\mu^2 \, ,\qquad\qquad v = \mu^{\frac1p}\phi \, ,
\end{equation}
factoring out the time dependence and recalling that we set $\beta =
\frac1{c_1}$, the kernel equation then looks
\begin{equation}
\label{e.K.r}
-\frac{a}{\mu^2}(\Delta \phi) + m\phi - |\phi|^{2p}\phi + R_V(\phi)=0 \, ,
\end{equation}
where 
\[
  R_V(\phi) = \frac1{\mu^{2+\frac1p}} \Pi_V 
  \quadr{N\tond{\mu^{\frac1p}\phi+w\tond{\mu^{\frac1p}\phi}} - N\tond{\mu^{\frac1p}\phi}}  \, .
\]

Since $R_V(\phi)$ is small, (see Lemma~\ref{l.small}), the kernel
equation, in the form \eqref{e.K.r}, appears as a perturbation of the
dNLS model studied in paper \cite{BamP09}. The main result of
\cite{BamP09} ensures the existence of breathers in the form of
non-degenerate ground states. We exploit non-degeneracy to continue such
solutions to solutions of the full equation \eqref{e.K.r}.

In order to exploit the result of \cite{BamP09} we need to work using
the norm $\Q$ rescaled by a factor $\mu^{\frac{n}2}$; so let us
introduce the following notations
\begin{equation}
\label{e.ell2r}
  \ell^2_\mu := \tond{\ell^2\,;\nlmu{\cdot}} \, ,
  \qquad\qquad \nlmu{\cdot}^2:=\mu^n\nll{\cdot}^2 \, ,
\end{equation}
and
\begin{equation}
\label{e.Qr}
  \Q_\mu := \tond{\ell^2\,;\nqmu{\cdot}} \, ,
  \qquad\qquad \nqmu{\cdot}^2:=\mu^n\nq{\cdot}^2 \, .
\end{equation}

In terms of these norms it is possible to prove good estimates for
$R_V$, see Lemma~\ref{l.small}. In particular we use them in order to
prove some discrete analogue of the Sobolev embedding theorems.

This ensure the applicability of the Implicit Function
Theorem~\ref{t.impl.K} for the continuation procedure (see
Proposition~\ref{t.K} for the details); we thus obtain the $2^n$
solutions $v^i=\mu^{\frac1p}\phi^i$, with $\phi^i$ close to $\psi^i$
(see estimate \eqref{est} and \eqref{e.est}), and the corresponding
$2^n$ solutions of \eqref{e.KG}
\[
  q^i(t) = v^i(\omega t) + w(v^i(\omega t)) \, , \qquad\qquad
  i=1,\dots,2^n \, .
\]

\subsubsection*{Step 3: Estimates} Finally one has to collect all the
estimates in order to get the result. We will also obtain the improved
estimate needed to control the sup norm of the difference between the
true solution and the actual solution. This will be done in
Sect.~\ref{s.eop} since it actually requires the results of all the
previous sections.  \qed


\section{The range equation}
\label{s.R}

In this Section we will prove Proposition~\ref{t.R} on the solutions
of the range equation. We start by controlling the inverse of the
linear operator $L^{(\omega)}$ defined in \eqref{e.LL}:
\begin{lemma}                   
\label{l.L}
If $0<a<\frac12$ and $|\omega^2-1|<
\frac12$, the linear operator $\left(L^{(\omega)}\right)^{-1}$ is
bounded from $W_0\subset X_0$ to $W_2\subset X_2$
\end{lemma}                     

\proof We use Neumann formula for the inversion of a linear
operator. We rewrite $L^{(\omega)}$ as
\[
L^{(\omega)} = L_\omega (\mathbb{I} + a L_\omega^{-1}\Delta) \, ,
\]
from which
\[
  \left(L^{(\omega)}\right)^{-1} = 
  (\mathbb{I} + a L_\omega^{-1}\Delta)^{-1} L_\omega^{-1} =
  \left[\sum_{k\geq 0}(-a L_\omega^{-1}\Delta)^k\right] L_\omega^{-1}
  \, .
\]
We observe that the series defines a bounded linear operator if
\begin{equation}
\label{e.Neu}
a\norma{L_\omega^{-1}\Delta}_{\L(W_0,W_2)}<1 \, .
\end{equation}
On one hand $\Delta$ acts only on the spatial index $j$ and defines a
bounded operator on $\ell^2$ with
\[
\nll{\Delta w}\leq 2\nll{w} \, ;
\]
on the other hand, $L_\omega^{-1}$ acts only on the temporal index $l$
and, provided $|\omega^2-1|<1/2$, is bounded. 
Hence \eqref{e.Neu} is fulfilled provided $a$ is small enough, i.e.
$a < \frac12$.

\qed                       

 We will show that $w$ is smaller than $v$ so that
 $N(v+w)=N(v)+$higher order terms. It is thus natural to expect
 the solution of the range equation \eqref{e.R} to be close to the
 solution of $L^{(\omega)} w = \Pi_W(N(v))$ namely to
\[
  w_0(v) := \tond{L^{(\omega)}}^{-1}\Pi_WN(v) \, .
\]

Define 
\begin{equation}
\label{udelta}
  \U_\delta:=
  \left\{v\in V_2\ :\  \nll{N(v)} + \norma{N'(v)}_{\L(\ell^2,X_0)} 
  <\delta  \right\}
\end{equation}
where we used the identification of $V_2$ with $\ell^2$ (see
\eqref{e.uvw}).

We are now ready to prove the following proposition

\begin{proposition}
\label{t.R}
Under the assumptions of Lemma \ref{l.L}, there exists $\delta>0$,
constants $C$, $C'$, and a function $w\in {\mathcal
  C}^{1,1}(\U_\delta,X_2)$ close to $w_0$ such that $w=w(v)$ solves
\eqref{e.R}. Moreover the following estimates hold
\begin{eqnarray}
\label{e.w}
\norma{w(v)-w_0(v)}_{X_2} &\leq& C \nll{N(v)}\norma{N'(v)}_{\L(\ell^2,X_0)}
\\
\label{e.w.1}
\norma{w'(v)}_{\L(X_2,X_2)} &\leq& C' \norma{N'(v)}_{\L(X_2,\ell^2)}.
\end{eqnarray}
\end{proposition}

\begin{remark}
Clearly $w(v)$ is small with $N(v)$, i.e.
\begin{equation}
\label{e.wnv}
\norma{w(v)}_{X_2}\leq C''\nll{N(v)} \; .
\end{equation}
\end{remark}

\proof Since the Nemitski operator defined by $N$ is ${\mathcal
  C}^{1,1}(X_2,X_0)$ by lemma \ref{l.smooth.N}, the implicit function
theorem ensures the existence of a neighbourhood of the origin in
which the function $w$ is well defined. In order to get the size of
such a neighbourhood and to prove the estimates \eqref{e.w},
\eqref{e.w.1} we go through the proof of the implicit function theorem
using the contraction mapping principle.

Let us first check that $w_0(v)\not\equiv 0$.  We write explicitly
$w_0(v)$:
\begin{eqnarray*}
 [ \Pi_WN(v)]_j &=& \Pi_W \left(|v_j|^{2p}v_j |e_1|^{2p}e_1\right) =
  |v_j|^{2p}v_j \frac1{2\pi} \sum_{l\in \L}
    \left[\int_0^{2\pi}|e_1|^{2p}e_1e_ldt \right] e_l = \\
  &=& |v_j|^{2p}v_j\frac1{2\pi}\sum_{l\in \L}c_l e_l \, ;
\end{eqnarray*}
%
%
%
since $c_l$ is the Fourier coefficients of $|\cos(t)|^{2p}\cos(t)$,
there surely exists at least one $c_l\neq 0$ with $l\in \L$.

Now rewrite the range equation as follows:
\begin{equation}
\label{e.R.1}
  w =\left(L^{(\omega)}\right)^{-1} \Pi_WN(v+w) =: F(v+w)\ . 
\end{equation}
Fix a positive $\delta_1$, and $v$ with $\|v\|<\delta_1$, and look for
conditions such that, the r.h.s. of \eqref{e.R.1} defines a
contraction of a ball of center $w_0(v)$ and radius $\delta_2$. We
claim that
\[
  \norma{w-w_0(v)}_{X_2} < \delta_2 \quad \Longrightarrow \quad
  \norma{F(v+w)-w_0(v)}_{X_2} < \delta_2
\]
if $\delta_2:= K \nll{N(v)}\norma{N'(v)}_{\L(\ell^2,X_0)}$ and $K$
sufficiently large.

To prove it, first remark that, since $N\in\C^{1,1}$, 
\begin{equation}
\label{e.R.11}
\sup_{\|u-v-w_0(v)\|_{X_2}\leq \delta_2}\| N'(u)\|_{\L(X_2,X_0)}\leq
\norma{N'(v)}_{\L(X_2,X_0)}+C'\delta_2
\, ,
\end{equation}
from which
\begin{eqnarray*}
\norma{F(v+w)-w_0(v)}_{X_2} &=& \norma{\left(L^{(\omega)}\right)^{-1}
\Pi_W\quadr{N(v+w)-N(v)} }_{X_2} \leq
\\
&\leq& C \sup_{\|u-v-w_0(v)\|_{X_2}\leq \delta_2} \|
N'(u)\|_{\L(X_2,X_0)} \|w\|_{X_2} \leq
\\
&\leq& C \tond{C'\delta_2+\norma{N'(v)}_{\L(X_2,X_0)}}
\tond{\delta_2+C''\nll{N(v)}} \, =
\\
&=& \frac{C}{K}\tond{C'K\nll{N(v)}+1}
\tond{K\norma{N'(v)}_{\L(X_2,X_0)}+C''} \delta_2 \, ,
\end{eqnarray*}
which is smaller than $\delta_2$ provided $K$ is sufficiently big and
$\nll{N(v)}$ and $\norma{N'(v)}_{\L(X_2,X_0)}$ sufficiently small,
i.e. if $v\in\U_\delta$. Then one immediately sees that by possibly
reducing $\delta$ the r.h.s. of \eqref{e.R.1} has a Lipschitz constant
smaller than one in the considered ball. So, we get the first of the
estimates \eqref{e.w}.

In order to get the estimate of the derivative of $w$ just remark that 
$$
w'(v)=\left(1-F'(v+w(v))\right)^{-1}F'(v+w(v))
$$
where we denoted $F(u):=\left(L^{(\omega)}\right)^{-1}
\Pi_WN(u) $. Using Neumann formula in order to compute
$\left(1-F'(v+w(v))\right)^{-1}$ one sees that this is a well defined
bounded linear operator provided $\delta$ is small enough. Adding the
estimate of $F'$ in the ball, which in turn is obtained through
\eqref{e.R.11} one gets the thesis. \qed

\begin{remark}
\label{c.contQ}
Since the topology induced by the $\Q$ norm is stronger than the
$\ell^2$ topology, one also has that the solution $w(v)$ of the range
equation is $\C^1(\Q\cap \U_\delta,X_2)$.
\end{remark}

In order to use in an effective way the inequality \eqref{e.wnv} we
will make use of the following lemma, which will be proved in Appendix
B, and which is based on the use of Sobolev embedding theorem applied
to functionals interpolating the $\ell^q$ norms.

\begin{lemma}
\label{nv}
One has 
\begin{equation}
\label{nv.0}
\mu^n \sum_{j\in\Z^n}|v_j|^{q} \leq C \nqmu{v}^{q} \; ,
\end{equation}
which gives in the case of operator $N$
\begin{equation}
\label{nv.1}
\norma{N(v)}_{\ell^2_\mu} \leq C'\nqmu{v}^{2p+1} \; .
\end{equation}
\end{lemma}

\begin{remark}
\label{nv.r}
The estimate \eqref{nv.1} is much stronger than the trivial one
obtained by using the homogeneity of $N$ and rescaling $v$. This will
be crucial for our development.
\end{remark}


\section{The kernel equation}
\label{s.Q}

We have seen in formula \eqref{e.K.r} that the kernel equation looks
like
\begin{equation}
\label{e.G}
G(\phi,\mu)=0,\qquad\qquad G(\phi,\mu) = G_0(\phi,\mu) + R_V(\phi,\mu)
\end{equation}
where
\begin{equation}
\label{e.G0}
G_0(\phi,\mu) := -\frac{a}{\mu^2}(\Delta \phi) + m\phi -
|\phi|^{2p}\phi=0
\end{equation}
is the equation for the ground state of the dNLS model. The maps $G$
and $G_0$ will be considered as maps $G:\Q_\mu\times \R\to \ell^2_\mu$.

The idea is to continue a solution of \eqref{e.G0} to a solution of
\eqref{e.G}. So, first we show that $R_V$ is actually a perturbation
of $G_0$. Denoting by $\Phi=\Phi(\mu)$ a solution of $G_0=0$, we then
show that $G_0'(\Phi)$ is an isomorphism of $\Q_\mu$ on $\ell^2_\mu$.

We begin by recalling the statement of the main result of
\cite{BamP09}. First we need to introduce a few objects.  Consider the
functional
\begin{equation}
\label{e.2}
H_0(\phi) := \mu^n\quadr{\frac12
  \sum_{|j-k|=1\,;\ j,k\in\Z^n}\frac{|\phi_j - \phi_k|^2}{\mu^{2}} -
  \frac1{p+1}\sum_{j\in\Z^n}|\phi_j|^{2p+2}},
\end{equation}
and the surface
\[
S=\left\{ \phi\ :\ N(\phi)=1 \right\}, \qquad N(\phi) :=
\mu^n\sum_{j\in\Z^n}|\phi_j|^2.
\]

\begin{theorem}
\label{t.BamP09} (Main theorem of \cite{BamP09}.)
For any $\mu$ small enough and $\frac12\leq p<2/n$ there exist $2^n$
distinct real valued sequences $\Phi_j^i(\mu)$ which are solutions of
\eqref{e.K.r} with $R_V\equiv 0$. Such solutions are coercive minima
of $H_0\big|_{S}$. The coercivity is intended in the norm
$\Q_\mu$. Furthermore one has
\begin{equation}
\label{est}
\nqmu{\Phi^i - \psi^i}\leq C\mu \, ,
\end{equation}
where $\psi_j^i$ are the sequences defined in \eqref{e.Psiref}.
\end{theorem}

\begin{remark}
In the following we will concentrate on one of these solutions, so we
will suppress the index $i$, from $\Phi$, from $\Psi$ and also from
$\psi$.
\end{remark}

\begin{lemma}
\label{l.small}
$R_V\in\C^{1,1}(\Q_\mu,\ell^2_\mu)$ fulfils
\begin{equation}
\label{e.small}
\nlmu{R_V(\phi)}\leq
C_1\mu^{2-\frac{n}2}\nqmu{\phi}^{4p+1} \, , \qquad\qquad 
\norma{R_V'(\phi)}_{\L(\Q_\mu,\ell^2_\mu)}\leq
C_2\sqrt\mu\nqmu{\phi}^{4p}.
\end{equation}
Moreover, $R_V(\phi,\mu)$ and $R'_V(\phi,\mu)$ are continuous with
respect to $\mu\in(0,+\infty)$.
\end{lemma}

\proof The smoothness of $R_V$ follows from the smoothness of $N$,
since both the norms in the spaces $\Q$ and $\ell^2$ have been
rescaled by the same factor.

We prove the first of \eqref{e.small} working first on the non
rescaled quantity
\begin{eqnarray*}
   \nll{\Pi_V\quadr{N(v+w(v)) -  N(v)}}^2 
   &\leq& \norma{N(v+w(v)) - N(v)}^2_{X_0}\leq
\\ &\leq& C\tond{\nll{\grad{N(v)}}^2 + \nll{\grad{N(w)}}^2}
   \norma{w(v)}^2_{X_2} \leq 
\\ &\leq& C_1\tond{\nll{\grad{N(v)}}^2 +
   \nll{\grad{N(w)}}^2 }\nll{N(v)}^2
\end{eqnarray*}
where we set $\grad{N(v)}:= |v|^{2p}$. We now apply the scaling
$v=\mu^{\frac1p}\phi$; and using Lemma~\ref{nv} one has
\begin{eqnarray*}
\nll{\grad{N(v)}}^2 &=& \sum_j|v_j|^{4p} \leq
\mu^{4-n} \nqmu{\phi}^{4p} \; ,
\\
\nll{\grad{N(w)}}^2 &=& \sum_j|w_j|^{4p} \leq C\nll{N(v)}^{4p}
\leq C_2 \mu^{4+\frac8p-2np} \nqmu{\phi}^{4p+8p^2} \; ,
\\
\nll{N(v)}^2 &=& \sum_j|v_j|^{4p+2} \leq
\mu^{4-n+\frac2p} \nqmu{\phi}^{4p+2} \; .
\end{eqnarray*}
If we ignore the term $\nll{\grad{N(w)}}$ which is much smaller in
$\mu$ than the main one $\nll{\grad{N(v)}}$ and we take the scaled
norm $\nlmu{\cdot}$ we obtain
\[
\nlmu{R_V(\phi)}\leq C_1\mu^{2-\frac{n}2}\nqmu{\phi}^{4p+1}.
\]
In order to conclude, we move to the estimate of $R'_V(v)$. Notice that
\[
\norma{R_V'(\phi)}_{\L(\Q_\mu,\ell^2_\mu)} =
\sup\frac{\nlmu{R_V'(\phi)[h]}}{\nqmu{h}}= \sup\frac{\nll{R_V'
    (\phi)[h]}}{\nq{h}} \leq \norma{R_V'(\phi)}_{\L(\Q,\ell^2)}.
\]

Let us first differentiate the non rescaled version of $R_V$:
\begin{equation}
\label{e.split.Q1}
  \frac{{\rm d}}{{\rm d} v} \Pi_V\circ[N(v+w(v)) - N(v)] = 
  \Pi_V\circ[N'(v+w(v)) - N'(v)] + \Pi_V\circ N'(v+w(v))\circ w'(v) \,,
\end{equation}

We deal with the first addendum in the r.h.s. of
\eqref{e.split.Q1}, by using the Lipschitz continuity of $N'$
\begin{eqnarray*}
\norma{\Pi_V\circ[N'(v+w(v)) - N'(v)]}_{\L(\Q,\ell^2)}&\leq&
\norma{[N'(v+w(v)) - N'(v)]}_{\L(X_2,X_0)}\leq \\ &\leq&
C(\norma{v}_{X_2}+\norma{w}_{X_2})^{2p-1}\norma{w}_{X_2}\leq \\ &\leq&
C_1(\nll{v}+\nll{N(v)})^{2p-1}\nll{N(v)} \, .
\end{eqnarray*}
Hence, after rescaling the variable and the norm, we get
\[
C\mu^4(\nll{\phi}+\mu^2\nll{N(\phi)})^{2p-1}\nll{N(\phi)}\leq C_1
\mu^2\mu^{2-np+\frac{n}2}\nqmu{\phi}^{4p}<C_1\mu^2\sqrt\mu\nqmu{\phi}^{4p}
\, ,
\]
where we have ignored the smaller term $\mu^2\nll{N(\phi)}$. A similar
estimate can be obtained for the second addendum in
\eqref{e.split.Q1}. Coming back to $R_V$
\[
  R_V'(\phi) = \mu^{-2} \frac{{\rm d}}{{\rm d} v} \Pi_V
  \quadr{N\tond{\mu^{\frac1p}\phi+w\tond{\mu^{\frac1p}\phi}} - 
    N\tond{\mu^{\frac1p}\phi}} \, ,
\]
the above informations yield
\[
\norma{R_V'(\phi)}_{\L(\Q_\mu,\ell^2_\mu)} \leq 
C_1\sqrt\mu\nqmu{\phi}^{4p} \, .
\]
The continuity with respect to the parameter $\mu$ follows from
standard arguments.
\qed

\begin{lemma}
\label{l.iso}
Let $\Phi $ be a coercive minimum of $H_0\big|_{S}$, then the
differential $G_0'(\Phi )$ is an isomorphism of $\Q_\mu$ on
$\ell^2_\mu$.
\end{lemma}

\proof By the theory of Lagrange multipliers one has that $\Phi$ is a
(free) critical point of $H_0+mN$, with a suitable $m$, while
$G'_0(\Phi)$ is such that 
$$
d^2(H_0+mN)(\Phi)[h,h]=\langle G'_0(\Phi)h;h\rangle_{\ell^2}\ .
$$
Introduce now ``polar coordinates'' 
$$
T_\Phi S\times \R\ni (\xi,\eta)\mapsto \phi(\xi,\eta):=\xi+\eta \Phi\ ,
$$
and write $G'_0(\Phi)$ as a block matrix in terms of such
coordinates. It has the structure
\begin{equation}
\label{block}
G'_0(\Phi)=
\left[\begin{matrix}
A&b
\\
b^T& d
\end{matrix}
\right]\ .\end{equation} By non-degeneracy one has $\langle
A\xi;\xi\rangle_{\ell^2}\geq C \nqmu{\xi}^2$, which, by Lax-Milgram
lemma implies that $A:\Q_\mu\to\ell^2_\mu$ is an isomorphism and is
positive definite.

Now one has 
\begin{equation}
\label{}
d=d^2(H_0+mN)(\Phi)[\Phi,\Phi]=\inter{G'(\Phi )[\Phi ],\Phi } =
-2p\sum_j{(\Phi )^{2p+2}} < 0
\end{equation}
which shows that the quadratic form $d^2(H_0+mN)(\Phi)[\Phi,\Phi]$ has
a negative direction. It follows that such a quadratic form does not
have null directions. Thus non-degeneracy and the thesis follow.
\qed

Lemmas \ref{l.small} and \ref{l.iso} give the hypothesis needed by
Theorem \ref{t.impl.K}. We get the final

\begin{proposition}
\label{t.K}
Let $\Phi(\mu)$ be a non-degenerate critical point of $H_0\big|_{S}$, then,
for $\mu$ small enough there exists a solution $\phi (\mu)$ of the
rescaled kernel equation \eqref{e.K.r}, such that
\begin{equation}
\label{e.est}
\nqmu{\phi - \Phi}\leq C\mu^{2-\frac{n}{2}}.
\end{equation}
\end{proposition}

\proof We verify the assumptions of Theorem \ref{t.impl.K}. Set
$X:=\Q_\mu$ and $Z:=\ell^2_\mu$ and $x_0(\mu):=\Phi (\mu)$. 
Define 
\[
F_0 := G_0(\phi,\mu),\qquad\qquad F_1 := R_V(\phi,\mu)\ .
\]
Then, by the previous lemmas the assumption of Theorem \ref{t.impl.K}
are satisfied and the thesis follows. \qed

\section{End of the proof}
\label{s.eop}

We begin a section by a simple lemma needed to obtain the estimate
\ref{e.supest}.

\begin{lemma}
\label{sup}
For any $j\in\Z^n$ we have
\begin{equation}
\label{grad.2}
\left|\phi_j\right|\leq 2\mu^{\frac12} \nq{\phi} \; .
\end{equation}
\end{lemma} 
\proof We write the proof for the case $n=2$. The case $n=1$ is
simpler. Denote $j=(h,k)$; one has
\begin{displaymath}
\phi_{h,k}^2 = \sum_{m=-\infty}^h(\phi_{m,k}^2 - \phi_{m-1,k}^2) =
\sum_{m=-\infty}^h(\phi_{m,k} - \phi_{m-1,k})(\phi_{m,k} + \phi_{m-1,k})
\end{displaymath}
which gives
\begin{eqnarray*}
\sup_{(h,k)\in\Z^2}\phi_{h,k}^2 &\leq& 4 \sqrt{\sum_{m\in\Z}\phi_{m,k}^2}
\sqrt{\sum_{m\in\Z}(\phi_{m+1,k} - \phi_{m,k})^2} \leq
\\
&\leq& \tond{2\sqrt{\norma
  {\phi}_{\ell^2}}\mu^{\frac{1}{2}-\frac{n}{4}} \left[ \mu^n
  \frac{\left\langle \phi;-\Delta\phi  \right\rangle_{\ell^2 }}{\mu^2}
  \right]^{1/4}}^2 
\end{eqnarray*}
and \eqref{grad.2}.\qed

\noindent {\it End of the proof of theorem \ref{t.main}.} We collect
now the estimates needed to conclude the proof. From proposition
\ref{t.K} and theorem \ref{t.BamP09} we get 
$$
\nqmu{\phi-\psi}\leq C\mu \; .
$$
From this we get (using also Lemma \ref{sup} and Theorem \ref{t.BamP09})
\begin{eqnarray*}
\nqmu{v}&\leq& C\mu^{\frac1p}
\\
\norma{v-\mu^{1/p}\psi}_{\ell^2}&\leq& \mu^{\frac1p-\frac{n}{2}+1}
\\
\sup_{j}\left|v_j-\mu^{1/p}\psi_j\right|&\leq&
C\mu^{\frac1p-\frac{n}{2}+\frac32} 
\\
\norma{w}_{X_2}&\leq& C\mu^{\frac{1}{p}-\frac{n}{2}+2}
\\
\left| w_j\right|&\leq& C\mu^{\frac{1}{p}-\frac{n}{2}+2}
\end{eqnarray*}
from which the thesis immediately follows. \qed


\appendix

\section{Smoothness of Nemitski operators.}
\label{s.nem}

\begin{lemma}
\label{l.nem.1}
If $p\geq \frac12$ then the operator $N$ defined in \eqref{e.NN}
is ${\C}^{1,1}$ from $\ell^2$ to $\ell^2$ with the usual norm.
\end{lemma}

\proof We first remark that $\ell^2\subset\ell^q$ for all $q\geq 2$,
indeed
\[
\norma{u}^q_{\ell^q} = \sum_j|u_j|^q = \sum_j(|u_j|^2)^{\frac{q}2}\leq
C_q(\sum_j{u_j^2})^\frac{q}2 = C_q\norma{u}_{\ell^2}^q,
\]
which immediately tells that $N(u)\in\ell^2$ if $u\in\ell^2$. Indeed
\[
\nll{N(u)}^2 = \sum_{j\in\Z^n}{|u_j|^{4p+2}} =
\norma{u}^{4p+2}_{\ell^{4p+2}}\leq C_p\nll{u}^{4p+2}.
\]
Moreover it will be useful to remind that $\ell^2\subset\ell^\infty$
since
\[
\norma{u}^2_{\ell^\infty} = (\sup_{j}|u_j|)^2 =
\max_j{|u_j|^2}\leq\sum_j{|u_j|^2} = \nll{u}^2.
\]
The first continuous embedding immediately gives also the continuity
at the origin. To obtain the continuity at a point $u\not=0$ we
proceed showing that $N$ is Frechet differentiable at any $u$ with
bounded differential $N'(u)$. From a direct computations one has
that
\[
N'(u)[h] = |u|^{2p} h,
\]
hence
\[
\nll{N'(u)[h]}^2 = \sum_j{|u_j|^{4p}h_j^2}\leq
C_1\nll{u}^{4p}\nll{h}^2
\]
which gives
\[
\norma{N'(u)}_{\L(\ell^2,\ell^2)}\leq C_p\nll{u}^{2p}.
\]
The possibility of locally bounding the differential yields to the
local Lipschitz continuity of $N$, since
\[
\nll{N(u)-N(v)}\leq
\tond{\sup_{\nll{w}\leq\nll{u}+\nll{v}}\norma{N'(w)}_{\L(\ell^2,\ell^2)}}
\nll{u-v}.
\]
Finally, the (local) Lipschitz continuity of $N'$ is due to $p\geq
\frac12$. Indeed
\[
\nll{N'(u)[h] - N'(v)[h]}^2 = \sum_j
\tond{|u_j|^{2p}-|v_j|^{2p}}^2h_j^2\leq C_p \nll{u-v}^2\nll{h}^2
\]
with $C_p=C(p,\norma{u}_{\ell^\infty},\norma{v}_{\ell^\infty})$ a
local constant.
\qed

\begin{corollary}
\label{c.nem}
If $p\geq \frac12$ then the operator $N$ is ${\C}^{1,1}(\Q,\ell^2)$.
\end{corollary}

We are now interested in regularity of $N$ as a map from $X_2$ to
$X_0$. We first state an auxiliary Lemma:

\begin{lemma}
\label{l.sm.norm}
Let us define $f(t):=\nll{u(t)}:[0,2\pi]\rightarrow W_0$. If $u\in
X_2$ then $f\in H^1([0,2\pi])$. More precisely one has
\begin{equation}
\label{e.sm.norm}
\norma{f}_{H^1}\leq \sqrt2\norma{u}_{X_2}.
\end{equation}
\end{lemma}

\proof Surely $f\in L^2(\I)$, where we set $\I:=[0,2\pi]$. By
differentiating $f^2(t)$ we get
\[
\dot f = \frac{1}{f}\tond{\inter{u,\dot u}},
\]
thus
\[
|\dot f(t)|^2\leq \frac{1}{\nll{u}^2}\tond{\nll{u}\nll{\dot u}}^2\leq
2\nll{\dot u}^2;
\]
hence
\[
\int_0^{2\pi}|\dot f(t)|^2dt\leq 2\int_0^{2\pi}\nll{\dot u}^2\leq
2\norma{u}_{X_2}^2
\]
which implies $f\in H^1(I)$. The estimate \eqref{e.sm.norm} follows
from
\[
\norma{f}_{H^1}^2 = \int_0^{2\pi}|f(t)|^2dt+\int_0^{2\pi}|\dot
f(t)|^2dt \leq 2\int_0^{2\pi}\quadr{\nll{u(t)}^2+\nll{\dot u(t)}^2}\leq
2\norma{u}^2_{X_2}.
\]

\begin{lemma}
\label{l.smooth.N}
The nonlinear operator $N$ is $\C^{1,1}(X_2,X_0)$.
\end{lemma}

\proof From Lemma \ref{l.nem.1} it follows that, for any fixed $t\in
I$ it holds
\[
\nll{N(u(t))}^2\leq C_1\nll{u(t)}^{4p+2},
\]
so from the Sobolev embedding $L^q(I)\hookrightarrow H^1(I),\,q\geq 2$
we get
\[
\int_0^{2\pi}\nll{N(u(t))}^2\leq C_1\int_0^{2\pi}\nll{u(t)}^{4p+2}\leq
C_2\norma{\nll{u}}^{4p+2}_{H^1} \leq C_2 \norma{u}_{X_2}^{4p+2},
\]
which simply gives
\[
\norma{N(u)}_{X_0}\leq C\norma{u}_{X_2}^{2p+1}.
\]
Let us consider now the Frechet differential
\[
N'(u):h\in X_2\mapsto N'(u)[h] = |u|^{2p}h\in X_0;
\]
once more, for any fixed $t\in I$ one has
\[
\nll{N'(u(t))[h(t)]}^2\leq C_1 \nll{u(t))}^{4p}\nll{h(t)}^2
\]
thus, again from Sobolev embeddings, we get
\[
\int_0^{2\pi}\nll{N'(u(t))[h(t)]}^2\leq C_1
\int_0^{2\pi}\nll{u(t))}^{4p}\nll{h(t)}^2\leq
C_2\norma{u}^{4p}_{X_2}\norma{h}^2_{X_2};
\]
we have so proved that the differential is locally bounded
\[
\norma{N'(u)}_{\L(X_2,X_0)}\leq C \norma{u}^{2p}_{X_2},
\]
and hence $N\in\C^1(0)$. Gathering the previous
results we deduce
\[
\norma{N(u)-N(v)}_{X_0}\leq
\tond{\sup_{w\in[u,v]}\norma{N'(w)}_{\L(X_2,X_0)}}\norma{u-v}_{X_2}
\]
which is the local Lipschitz continuity. The local Lipschitz
continuity of $N'(u)$ can be obtained in the same way as in Lemma
\ref{l.nem.1}. Indeed
\begin{eqnarray*}
\norma{(N'(v+w)-N'(v))[h]}_{X_0}^2 &=&
\int_0^{2\pi}\nll{h(|v+w|^{2p}-|v|^{2p})}^2
\leq\\ &\leq&\int_0^{2\pi}\nll{h}^2\nll{|v+w|^{2p}-|v|^{2p}}^2\leq\\ &\leq&
\norma{\nll{h}}_{L^4}\norma{\nll{|v+w|^{2p}-|v|^{2p}}^2}_{L^2}\leq\\ &\leq&
C\norma{h}_{X_2}^2\norma{\nll{|v+w|^{2p}-|v|^{2p}}^2}_{L^2}.
\end{eqnarray*}
Moreover, following Lemma \ref{l.nem.1} one has
\[
\nll{|v+w|^{2p}-|v|^{2p}}^2\leq C_L(t)^2\nll{w}^2
\]
with
\begin{eqnarray*}
C_L(t) &=&
2p(\norma{v(t)}_{\ell^\infty}+\norma{w(t)}_{\ell^\infty})^{2p-1}\leq\\ &\leq&
4(\nll{v(t)}+\nll{w(t)})^{2p-1}\leq \\ &\leq&
4(\norma{\nll{v}}_{L^\infty}+\norma{\nll{w}}_{L^\infty})^{2p-1}\leq\\ 
&\leq& 4(\norma{v}_{X_2}+\norma{w}_{X_2})^{2p-1}.
\end{eqnarray*}
So
\[
\norma{N'(v+w)-N'(v)}_{\L(X_2,X_0)}\leq L\norma{w}_{X_2},\qquad\qquad
L:=C(\norma{v}_{X_2}+\norma{w}_{X_2})^{2p-1}
\]
where the Lipschitz constant $L$ is local, unless for $p=\frac12$.
\qed

\bigskip
\bigskip

\section{Approximation of discrete functionals and proof of Lemma
  \ref{nv}.}
\label{s.A.norms}

The proof of Lemma \ref{nv} is based on the use of Sobolev embedding
theorem applied to the continuous interpolation of some discrete
functional. In turn, following \cite{BamP09}, the continuous
interpolation is obtained through the method of finite elements as we
are now going to recall.

\subsection{The case $n=1$}

Define the sequence of functions $s_j(x)$ by
\begin{equation}
\label{fe.1}
s_j(x)=
\begin{cases}
0,\qquad\qquad  {\rm if}\quad |x-j|>1\\
x-j+1,\quad  {\rm if}\quad   -1\leq x- j\leq 0\\
- x+j+1,\quad  {\rm if}\quad  0\leq x-j\leq 1
\end{cases}
\end{equation}
and, to a sequence $\psi_j$, we associate a function
\begin{equation}
\label{e.5}
\Upsilon(x) := \sum_j\psi_js_j(x/\mu).
\end{equation}
On the interval $T_j:=[\mu j,\mu(j+1))$ the above function reads
\begin{equation}
\label{e.interpol}
\Upsilon(x) = (x-\mu j)\frac{(\psi_{j+1} - \psi_{j})}{\mu} + \psi_{j}.
\end{equation}

\subsection{The case $n=2$.}

For each multi-index $j=(h,k)$, take the function $s_{h,k}(x,y)$ which
represents the hexagonal pyramid of height one centered in $(h,k)$
with support the union of the six triangles
$T^{\pm}_{h,k},\,T^+_{h-1,k},\,T^-_{h,k+1},\,T^+_{h,k-1},\,T^-_{h+1,k}$. For
example, on $T^+_{h,k}$ the function $s_{h,k}$ represents the plane in
$\mathbb{R}^3$ which passes through the three points
$(h,k,1),\,(h+1,k,0),\,(h,k+1,0)$, namely
\begin{displaymath}
s_{h,k}(x,y) = -{x} -{y}+h+k+1.
\end{displaymath}
We take $\{s_{h,k}(x/\mu,y/\mu)\}_{(h,k)\in\Z^2}$ as a basis to
generate a piecewise linear function $\psi(x,y)$ which interpolates
$\psi_{h,k}$
\begin{equation}
\label{e.int.2}
\Upsilon(x,y) := \sum_{(h,k)\in\Z^2}{\psi_{h,k}s_{h,k}(x/\mu,y/\mu)}.
\end{equation}
On the triangle $T^+_{h,k}$ the function $\psi$ is a plane which reads
\begin{equation}
\Upsilon(x,y) = (x-\mu h)\frac{(\psi_{h+1,k} - \psi_{h,k})}\mu + (y-\mu
k)\frac{(\psi_{h,k+1} - \psi_{h,k})}\mu + \psi_{h,k},
\end{equation}
while on the opposite triangle $T^-_{h,k}$ it reads
\begin{equation}
\Upsilon(x,y) = (x-\mu h)\frac{(\psi_{h,k} - \psi_{h-1,k})}\mu + (y-\mu
k)\frac{(\psi_{h,k} - \psi_{h,k-1})}\mu + \psi_{h,k}.
\end{equation}

\begin{definition}
\label{linspace.2}
We denote by $\E_\mu$ the linear subspace of $H^1$ of the functions
\eqref{e.int.2} with $\{\psi_j\}\in\Q_\mu$.
\end{definition}

\subsection{Interpolation}

We recall now some lemmas which were proved in \cite{BamP09}.

\begin{lemma}
\label{grad}
Let $\Upsilon\in \E_\mu$, denote by $\psi=\{\psi_j\}$ the corresponding
sequence, then one has 
\begin{equation}
\label{grad.1}
\int_{\R^n}\left| \nabla
\Upsilon(x)\right|d^nx=-\frac{1}{\mu^2}\mu^n\langle
\psi,\Delta\psi\rangle_{\ell^2}
\end{equation}
\end{lemma}

\begin{lemma}
\label{l.pot.2d}
Let $\Upsilon\in \E_\mu$, denote by $\psi=\{\psi_j\}$ the corresponding
sequence; define
\begin{displaymath}
G_c(\Upsilon):=\int_{\R^2}|\Upsilon|^{q+2},\qquad
G_d(\Upsilon):=\mu^n\sum_{j\in\Z^n}{|\psi_j|^{q+2}},\qquad
R_G(\Upsilon):= G_c(\Upsilon) - G_d(\Upsilon);
\end{displaymath}
if $q\geq 1$ then $R_G\in\C^2(\E_\mu)$ and for any bounded open set
$\U\subset\E_\mu$ there exists $C(\U)$ such that
\begin{displaymath}
\norma{R_G}_{\C^2(\U)}\leq C \mu.
\end{displaymath}
\end{lemma}

\noindent {\it Proof of lemma \ref{nv}}. One has 
\begin{eqnarray*}
\mu^n\sum_{j}\left|\psi_j\right|^{q+2}=\int_{\R^n}\left|\Upsilon\right|^{q+2}
d^nx +O(\mu) \leq C\norma{\Upsilon}_{H^1}^{q+2}+O(\mu) 
\\
= C\nqmu{\psi}^{q+2}+O(\mu)\leq C'\nqmu{\psi}^{q+2} 
\end{eqnarray*} 
where the third inequality follows from the fact that, for $n=1,2$ one
has that the injection $H^1\hookrightarrow L^p$ is continuous for any
$p$. From this the thesis immediately follows. \qed

\section{A version of the implicit function theorem.}
\label{app.3}

\begin{theorem}
\label{t.impl.K}
Let $X,Z$ be Banach spaces and let $0\in\I\subset\R $. Let $F_0\in {\mathcal
  C}^{1,1}(\U\times\I,Z)$, with $\U\subset X$ open. 
Let $x_0=x_0(\mu)\in\U$ be such that 
\begin{enumerate}
\item 
\begin{equation}
\label{e.appr}
F_0(x_0(\mu),\mu)=0; \forall \mu\in\I
\end{equation}

\item $F'_{0,x}(x,\mu)$ is Lipschitz in $\U$ uniformly in $\mu$,
  i.e. there exists $L=L(\U)$ {\bf independent of $\mu$} such that
\begin{equation}
\label{e.Lip.D}
\norma{F'_{0,x}(x,\mu) - F'_{0,x}(x_0,\mu)}_{\L(X,Z)}\leq
L\norma{x-x_0}_X,\qquad\forall x\in\U;
\end{equation}

\item $F'_{0,x}(x_0,\mu)$ is invertible and its inverse is bounded
  uniformly in $\mu$, i.e. there exists $C_1(\U)$ such that
  \begin{equation}
\label{e.invbound}
\norma{[F'_{0,x}(x_0,\mu)]^{-1}}_{\L(Z,X)}<C_1 \,,\qquad\forall \mu\in\I
\, .
\end{equation}

\item Let $F_1\in {\mathcal C}^{1}(\U\times\I,Z)$ be such that there exist
  $\alpha>0$ and $C_2(\U)$ such that
\begin{equation}
\label{e.sm.F1}
\norma{F_1(x,\mu)}_Z \leq C_2\mu^\alpha,\qquad\forall
(x,\mu)\in\U\times\I;
\end{equation}

\item there exist $\beta>0$ and $C_3(\U)$ such that
\begin{equation}
\label{e.sm.dF1}
\norma{F'_{1,x}(x,\mu)}_{\L(X,Z)}<C_3 \mu^\beta,\qquad\forall
(x,\mu)\in\U\times\I;
\end{equation}

\end{enumerate}
Define $F:=F_0+F_1$, then there exist $\U_\delta\subset\U$ and
$\mu^*(\delta)<\mu$ and a function
$x(\mu):\I_0^*:=(0,\mu^*)\mapsto\U_\delta$ which solves
\begin{equation}
\label{e.fixedpoint}
F(x(\mu),\mu)=0 \, ,
\end{equation}
with
\begin{equation}
\label{e.sfia}
  \norma{x(\mu)-x_0(\mu)}_X\leq C\mu^\alpha \, .
\end{equation}
Moreover one has $x(\mu)\in\C^0(\I_0^*,\U_\delta)$.
\end{theorem}

\proof First remark that, provided $\mu$ is small enough and possibly
restricting $\U$,
$[F'_x(x)]^{-1}$ exists and fulfills 
\begin{equation}
\label{e.invbound.2}
\norma{[F'_{x}(x_0,\mu)]^{-1}}_{\L(Z,X)}<C \,,\qquad\forall
\mu\in\I\ ,
\, .
\end{equation}
Define now $A(x,\mu):=x- [F'_{x}(x_0,\mu)]^{-1}F(x,\mu)$ and remark
that any fixed point of $A$ is a solution of our problem. We now prove
that $A$ is a contraction of a ball of radius $O(\mu^\alpha)$ and
center $x_0$. So, let $x$ be such that $\norma{x-x_0}\leq \delta$, and
let us estimate the Lipschitz constant of $A$ in such a ball. One has
\begin{eqnarray*}
\norma{F'_{0,x}(x_0,\mu)-F'_{0,x}(x,\mu) }\leq C\delta
\\
\norma{F'_{1,x}(x_0,\mu)-F'_{1,x}(x,\mu) }\leq C\mu^\beta
\end{eqnarray*}
It follows
that 
$$
\norma{A'(x,\mu)}=\norma{[F'_x(x_0,\mu)]^{-1}(F'_x(x_0,\mu)-F'_x(x,\mu)
  )}\leq q<1
$$ 
provided $\mu$ and $\delta$ are small enough. Compute now 
\begin{eqnarray*}
\norma{A(x,\mu)-x_0}\leq
\norma{A(x_0,\mu)-x_0}+\norma{A(x,\mu)-A(x_0,\mu)} 
\\
\leq C\mu^\alpha q\delta\ ,
\end{eqnarray*}
which is smaller than $\delta$ provided $C\mu^\alpha<(1-q)\delta$,
which in turn can be obtained e.g. by taking 
$$
\delta=\frac{C\mu^\alpha}{2(1-q)}\ .
$$
  \qed


\section*{Acknowledgements}

Partially supported by PRIN 2007B3RBEY ``Dynamical Systems and applications''.


\begin{thebibliography}{10}

\bibitem{AubKK01}
{\sc S.~Aubry, G.~Kopidakis, and V.~Kadelburg}, {\em Variational proof for hard
  discrete breathers in some classes of {H}amiltonian dynamical systems},
  Discrete Contin. Dyn. Syst. Ser. B, 1 (2001), pp.~271--298.

\bibitem{Bak04}
{\sc S.~N. Bak}, {\em The constrained minimization method in the problem of the
  oscillations of a chain of nonlinear oscillators}, Mat. Fiz. Anal. Geom., 11
  (2004), pp.~263--273.

\bibitem{BakP04}
{\sc S.~N. Bak and A.~A. Pankov}, {\em On periodic oscillations of an infinite
  chain of linearly coupled nonlinear oscillators}, Dopov. Nats. Akad. Nauk
  Ukr. Mat. Prirodozn. Tekh. Nauki,  (2004), pp.~13--16.

\bibitem{Bam96}
{\sc D.~Bambusi}, {\em Exponential stability of breathers in {H}amiltonian
  networks of weakly coupled oscillators}, Nonlinearity, 9 (1996),
  pp.~433--457.

\bibitem{BamCP09}
{\sc D.~Bambusi, A.~Carati, and T.~Penati}, {\em Boundary effects on the
  dynamics of chains of coupled oscillators}, Nonlinearity, 22 (2009),
  pp.~923--946.

\bibitem{BamCP02}
{\sc D.~Bambusi, A.~Carati, and A.~Ponno}, {\em The nonlinear {S}chr\"odinger
  equation as a resonant normal form}, Discrete Contin. Dyn. Syst. Ser. B, 2
  (2002), pp.~109--128.

\bibitem{BP}
{\sc D.~Bambusi and S.~Paleari}, {\em Families of periodic solutions of
  resonant {P}{D}{E}s}, J. Nonlinear Sci., 11 (2001), pp.~69--87.

\bibitem{BP2}
\leavevmode\vrule height 2pt depth -1.6pt width 23pt, {\em Families of periodic
  orbits for some {PDE}'s in higher dimensions}, Commun. Pure Appl. Anal., 1
  (2002), pp.~269--279.

\bibitem{BamP09}
{\sc D.~Bambusi and T.~Penati}, {\em Continuous approximation of ground states
  in {DNLS} lattices},  (2009).
\newblock preprint.

\bibitem{BPon05}
{\sc D.~Bambusi and A.~Ponno}, {\em On metastability in {F}{P}{U}}, Comm. Math.
  Phys., 264 (2006), pp.~539--561.

\bibitem{BerL83}
{\sc H.~Berestycki and P.-L. Lions}, {\em Nonlinear scalar field equations.
  {I}. {E}xistence of a ground state}, Arch. Rational Mech. Anal., 82 (1983),
  pp.~313--345.

\bibitem{ColGM78}
{\sc S.~Coleman, V.~Glaser, and A.~Martin}, {\em Action minima among solutions
  to a class of {E}uclidean scalar field equations}, Comm. Math. Phys., 58
  (1978), pp.~211--221.

\bibitem{Dui84}
{\sc J.~J. Duistermaat}, {\em Bifurcation of periodic solutions near
  equilibrium points of {H}amiltonian systems}, in Bifurcation theory and
  applications (Montecatini, 1983), Springer, Berlin, 1984, pp.~57--105.

\bibitem{FriP99}
{\sc G.~Friesecke and R.~L. Pego}, {\em Solitary waves on {FPU} lattices. {I}.
  {Q}ualitative properties, renormalization and continuum limit}, Nonlinearity,
  12 (1999), pp.~1601--1627.

\bibitem{FriP1-4}
\leavevmode\vrule height 2pt depth -1.6pt width 23pt, {\em Solitary waves on
  {F}ermi-{P}asta-{U}lam lattices. {I}. {Q}ualitative properties,
  renormalization and continuum limit. {II}. {L}inear implies nonlinear
  stability. {III}. {H}owland-type {F}loquet theory. {IV}. {P}roof of stability
  at low energy}, Nonlinearity, 12/15/17 (1999/2002/2004),
  pp.~1601--1627/1343--1359/207--227/229--251.

\bibitem{GriSS87}
{\sc M.~Grillakis, J.~Shatah, and W.~Strauss}, {\em Stability theory of
  solitary waves in the presence of symmetry. {I}}, J. Funct. Anal., 74 (1987),
  pp.~160--197.

\bibitem{HofW08}
{\sc A.~Hoffman and C.~E. Wayne}, {\em Counter-propagating two-soliton
  solutions in the {F}ermi-{P}asta-{U}lam lattice}, Nonlinearity, 21 (2008),
  pp.~2911--2947.

\bibitem{Ioo00}
{\sc G.~Iooss}, {\em Travelling waves in the {F}ermi-{P}asta-{U}lam lattice},
  Nonlinearity, 13 (2000), pp.~849--866.

\bibitem{IooJ05}
{\sc G.~Iooss and G.~James}, {\em Localized waves in nonlinear oscillator
  chains}, Chaos, 15 (2005), pp.~015113, 15.

\bibitem{IooK00}
{\sc G.~Iooss and K.~Kirchg{\"a}ssner}, {\em Travelling waves in a chain of
  coupled nonlinear oscillators}, Comm. Math. Phys., 211 (2000), pp.~439--464.

\bibitem{IooP06}
{\sc G.~Iooss and D.~E. Pelinovsky}, {\em Normal form for travelling kinks in
  discrete {K}lein-{G}ordon lattices}, Phys. D, 216 (2006), pp.~327--345.

\bibitem{Jam03}
{\sc G.~James}, {\em Centre manifold reduction for quasilinear discrete
  systems}, J. Nonlinear Sci., 13 (2003), pp.~27--63.

\bibitem{JamS05}
{\sc G.~James and Y.~Sire}, {\em Travelling breathers with exponentially small
  tails in a chain of nonlinear oscillators}, Comm. Math. Phys., 257 (2005),
  pp.~51--85.

\bibitem{JamS08}
\leavevmode\vrule height 2pt depth -1.6pt width 23pt, {\em Center manifold
  theory in the context of infinite one-dimensional lattices}, in The
  {F}ermi-{P}asta-{U}lam problem, vol.~728 of Lecture Notes in Phys., Springer,
  Berlin, 2008, pp.~208--238.

\bibitem{Kal89}
{\sc L.~A. Kalyakin}, {\em Long-wave asymptotics. {I}ntegrable equations as the
  asymptotic limit of nonlinear systems}, Uspekhi Mat. Nauk, 44 (1989),
  pp.~5--34, 247.

\bibitem{KirSM92}
{\sc P.~Kirrmann, G.~Schneider, and A.~Mielke}, {\em The validity of modulation
  equations for extended systems with cubic nonlinearities}, Proc. Roy. Soc.
  Edinburgh Sect. A, 122 (1992), pp.~85--91.

\bibitem{MacA94}
{\sc R.~S. MacKay and S.~Aubry}, {\em Proof of existence of breathers for
  time-reversible or {H}amiltonian networks of weakly coupled oscillators},
  Nonlinearity, 7 (1994), pp.~1623--1643.

\bibitem{MW89}
{\sc J.~Mawhin and M.~Willem}, {\em Critical point theory and {H}amiltonian
  systems}, Springer-Verlag, New York, 1989.

\bibitem{MizP08}
{\sc T.~Mizumachi and R.~L. Pego}, {\em Asymptotic stability of {T}oda lattice
  solitons}, Nonlinearity, 21 (2008), pp.~2099--2111.

\bibitem{SchW00}
{\sc G.~Schneider and C.~E. Wayne}, {\em Counter-propagating waves on fluid
  surfaces and the continuum limit of the {F}ermi-{P}asta-{U}lam model}, in
  International Conference on Differential Equations, Vol. 1, 2 (Berlin, 1999),
  World Sci. Publishing, River Edge, NJ, 2000, pp.~390--404.

\bibitem{SepM97}
{\sc J.-A. Sepulchre and R.~S. MacKay}, {\em Localized oscillations in
  conservative or dissipative networks of weakly coupled autonomous
  oscillators}, Nonlinearity, 10 (1997), pp.~679--713.

\bibitem{Sir05}
{\sc Y.~Sire}, {\em Travelling breathers in {K}lein-{G}ordon lattices as
  homoclinic orbits to {$p$}-tori}, J. Dynam. Differential Equations, 17
  (2005), pp.~779--823.

\bibitem{Wei99}
{\sc M.~I. Weinstein}, {\em Excitation thresholds for nonlinear localized modes
  on lattices}, Nonlinearity, 12 (1999), pp.~673--691.

\end{thebibliography}

\def\cprime{$'$}
\def\i{\ii}

\end{document}